\documentclass[authoryear,hyperref,preprint,1p]{elsarticle}
\input{Fiskdef2.tex}
\ifx\StProb\undefined
\makeatletter
\def\ps@pprintTitle{%
     \let\@oddhead\@empty
     \let\@evenhead\@empty
     \def\@oddfoot{}%
     \let\@evenfoot\@oddfoot}
\makeatother
\fi
\begin{document}

\begin{frontmatter}
\author[tukl,itwm]{Peter Ruckdeschel}
\ead{Peter.Ruckdeschel@itwm.fraunhofer.de}
\author[ubt]{Helmut Rieder}
\ead{Helmut.Rieder@uni-bayreuth.de}
\address[tukl]{Fraunhofer ITWM, Abt.\ Finanzmathematik,
         Fraunhofer-Platz 1, 67663 Kaiserslautern, Germany}
\address[itwm]{TU Kaiserslautern, AG Statistik, FB.\ Mathematik,
         P.O.Box 3049, 67653 Kaiserslautern, Germany}
\address[ubt]{Mathematical Institute, University of Bayreuth,
              95440 Bayreuth, Germany
              }
\date{\today}

\title{Fisher Information of Scale}

\begin{abstract} 
Motivated by the information bound for the asymptotic variance of M-estimates for scale,
we define Fisher information of scale of any distribution function~$F$ on the real line
as the supremum of all
$  
   \bigl( \int x \:\phi'(x)\,F(dx)\bigr)^2
                          \big/\!\int \phi^2(x)\,F(dx)  $,
where $\phi$ ranges over the continuously differentiable functions with derivative of
compact support and where, by convention, $0/0:=0$. In addition, we enforce equivariance
by a scale factor. Fisher information of scale is weakly lower semicontinuous and convex.
It is finite iff the usual assumptions on densities hold, under which Fisher information
of scale is classically defined, and then both classical and our notions agree.
Fisher information of scale finite is also equivalent to $L_2$-differentiability and
local asymptotic normality, respectively, of the scale model induced by~$F$.
\end{abstract}
\end{frontmatter}
\section{Motivation and Definition}\label{setup}
If $F$ is any distribution function on~$\R$, the real line, and $\phi \colon \R \to \R$
a suitable scores function such that $\int \phi \,dF=0$, an M-estimate of scale $S_n$
may formally be defined by
\begin{equation}\label{Msc}
    \sum_{i=1}^n \phi \Bigl(\frac{x_i}{S_n}\Bigr) =0\;.
\end{equation}
The estimand refers to the scale model $(F_\sigma)_{0<\sigma<\infty}$
induced by~$F=F_1$, where $F_\sigma(x)=F(x/\!\sigma)$.
\par
Taylor expanding
  $ \phi(x/s)=\phi(x/\!\sigma) - (s-\sigma)\phi'(x/\!\sigma)\,x/\!\sigma^2 + \cdots $,
  we formally obtain
\begin{equation}
    \sqrt{n}\,(S_n-\sigma)= \sigma\frac{n^{-1/2} \sum_1^n \phi(x_i/\!\sigma)}{
                           n^{-1} \sum_1^n \phi'(x_i/\!\sigma)\,x_i/\!\sigma
    } \: + \: \cdots
\end{equation}
such that under observations $x_1, \dots,x_n $ i.i.d.$\sim F_\sigma$ and assuming
sufficient regularity, in particular consistency, $\sqrt{n}\,(S_n-\sigma)$ will as $n\to \infty$
be asymptotically normal with mean zero and variance 
\begin{equation}\label{asyVarMsc}
      V(\phi,F_\sigma)=\sigma^2 \;V_1(\phi,F)\,,\quad
      V_1(\phi,F):=
      \frac{     \int \phi^2(x)\,F(dx)  }{
                      \bigl( \int x \:\phi'(x)\,F(dx)\bigr)^2  }\;.
\end{equation}
If $\phi$ is differentiable with continuous derivative of compact support, both $\phi(x)$
and $x\,\phi'(x)$ are bounded, so the integrals in~\eqref{asyVarMsc} are well-defined
for any distribution~$F$ on the Borel $\sigma$-algebra~$\B$ of~$\R$.
As in the theory of generalized functions (\citet[Ch.~6]{Ru:91}), regularity conditions
are shifted to the test functions whenever possible.
\par
The usual information bound for asymptotic variance would say that
  $V(\phi,F_\sigma)\ge {\cal I}_{\rm s}^{-1}(F_\sigma)$ and, hopefully, the lower bound will
also be achieved.
\par
This leads us to the following definition of~${\cal I}_{\rm s1}(F)$.
The extension to~${\cal I}_{\rm s}(F_\sigma)$ for the scale transforms~$F_\sigma$
of~$F$ matches~\eqref{asyVarMsc}.
\begin{Def} \label{FiSc}
Fisher information of scale, for any distribution~$F$ on the real line, is defined by
\begin{equation}
  {\cal I}_{\rm s1}(F):= \sup_{\phi\in {\cal C}_{\rm c1}}
  \frac{ \bigl( \int x \:\phi'(x)\,F(dx)\bigr)^2}{\int \phi^2(x)\,F(dx)  }\;,
\end{equation}
where ${\cal C}_{\rm c1}$ denotes the set of all differentiable functions $\phi \colon\R\to\R$
whose derivative is continuous and of compact support, and $0/0:=0$ by convention.
For the scale transforms $F_\sigma$ of~$F$ we define
\begin{equation} \label{Isigdef}
   {\cal I}_{\rm s}(F_\sigma):= \sigma^{-2}{\cal I}_{\rm s1}(F)\,,\qquad 0<\sigma<\infty\,.
\end{equation}
\end{Def}
\begin{Rem}\rm
Since the map $\phi \mapsto \phi_\sigma$, where $\phi_\sigma(x):=\phi(\sigma x)$ and
  $\phi'_\sigma(x)=\sigma \,\phi'(\sigma x)$, defines a one-to-one correspondence
  on~${\cal C}_{\rm c1}$, we obtain scale invariance of~${\cal I}_{\rm s1}$,
\begin{equation}\label{Is1invariant}
  {\cal I}_{\rm s1}(F_\sigma)={\cal I}_{\rm s1}(F)\,,\qquad 0<\sigma<\infty\,.
\end{equation}
So extension~\eqref{Isigdef} is needed to obtain scale equivariance. In the scale model,
as opposed to location, it matters whether a given distribution~$F$ is considered
element~$F=F_1$ or, for example, element $F=F_{.5}$ (in the scale model
generated by~$F_2$). \qed\end{Rem}
\par
Motivated by the information bound, Definition~\ref{FiSc} is instrinsically statistical.
It does not a priori use the assumption of, and suitable conditions on, densities.
These properties rather follow from the definition in case ${\cal I}_{\rm s}$
is finite. Another advantage is that Definition~\ref{FiSc} implies certain
topological properties (convexity and lower continuity) of~${\cal I}_{\rm s}$.
\par
The definition parallels \citet[Def.~4.1]{Hu:81} in the location case,
\begin{equation}\label{HuberIFloc}
  {\cal I}_{\rm l}(F):= \sup_\phi  \frac{ \bigl( \int \phi'(x)\,F(dx)\bigr)^2
                          }{     \int \phi^2(x)\,F(dx)  }\;,
\end{equation}
where $\phi $, subject to $\int \phi^2 dF>0$, ranges over the (smaller) set
 ${\cal C}_{\rm c}^1$ of all continuously differentiable functions which
 themselves are of compact support. ${\cal I}_{\rm l}$ is shift invariant.
\par
\citet[p.~79]{Hu:81}, states vague lower semicontinuity and convexity
of~${\cal I}_{\rm l}$. By \citet[Thm.~4.2]{Hu:81}, ${\cal I}_{\rm l}(F)$
is finite iff $F$ is absolutely continuous with an absolutely continuous
density~$f$ such that $f'\!\!/\!f \in L_2(F)$, in which case
   ${\cal I}_{\rm l}(F)=\int (f'\!\!/\!f)^2\,dF$.
\begin{Rem}\rm
The latter result, by arguments of the proof to Theorem 2.2 below, still obtains if
definition~\eqref{HuberIFloc} is based on ${\cal C}_{\rm c1}$. Only vague lower
semicontinuity of~${\cal I}_{\rm l}$ would be weakened to weak continuity
(which, however,  makes no difference in the setup of normed measures).
The convention $0/0:=0$ could replace the side condition $\phi\ne0$ a.e.~$F$
in~\eqref{HuberIFloc} as well.
\par
The non-suitability of~${\cal C}_{c}^1$, and suitability of ${\cal C}_{\rm c1}$
instead, is the tribute to the scale model, for which the functions
  $x \mapsto x\,\phi'(x)$ need to be dense in $L_1(F_0)$ with respect to the
punctuated (substochastic) measure~$F_0$ introduced in~\eqref{F0} below.
\qed\end{Rem}
\par
Fisher information of scale has been treated by \citet{Hu:64,Hu:81} not in the
previous generality but only under suitable assumptions on densities and,
in an auxiliary way, has been reduced to the location case by symmetrization
and the log-transform, \citet[Sec.~5.6]{Hu:81}.
\section{Main Results}
\begin{Prop}\label{vsc+cvx}
    ${\cal I}_{\rm s1}$ is weakly lower semicontinuous and convex.
\end{Prop}
\par Zero observations do not contain any information about scale.
Removing the mass of any distribution~$F$ at zero, we define the punctuated,
possibly substochastic measure~$F_0$ by
\begin{equation}\label{F0}
   F_0:=F - F(\{0\})1_0\;,
\end{equation}
where~$1_0$ denotes Dirac measure at~$0$. In terms of distribution functions,
denoting by $1_{[0,\infty)}$ the indicator function, we have
  $F_0(x)=F(x) - (F(0)- F(0-))1_{[0,\infty)}(x)$.
\begin{Thm}\label{Thmscal} 
For any distribution~$F$ on the real line, ${\cal I}_{\rm s1}(F)$ is finite iff
\begin{enumerate}
\item[i)] $F_0$ is absolutely continuous with a density~$f$ such that
\item[ii)] $x \mapsto x\,f(x)$ is absolutely continuous, and
\item[iii)] $x \mapsto \Lambda(x):=-[\,xf(x)]'\!/\!f(x) \in L_2(F_0)$,
\end{enumerate}
in which case
    $ {\cal I}_{\rm s1}(F)= \int  \Lambda^2 \,dF_0 =
      {\displaystyle\int_{x\ne0}} [\, 1 + x f'(x)/\!\!f(x)\,]^2\,F(dx)
    $\,.
\end{Thm}
\section{Consequences for the Scale Model}\label{conseq}
\noindent For the scale transforms~$F_\sigma$ of~$F$,
  $ {\cal I}_{\rm s1}(F_\sigma)={\cal I}_{\rm s1}(F) $ and
  ${\cal I}_{\rm s}(F_\sigma)= \sigma^{-2}{\cal I}_{\rm s1}(F)$ by
  \eqref{Is1invariant} and~\eqref{Isigdef}, respectively.
In particular, ${\cal I}_{\rm s1}(F_\sigma)$ and ${\cal I}_{\rm s}(F_\sigma)$
are finite iff ${\cal I}_{\rm s1}(F)$ is finite.
Also conditions i) and~ii) of Theorem~\ref{Thmscal} are simultaneously fulfilled
for a density~$f$ of~$F_0$ and the density $f_\sigma(x)=\sigma^{-1}f(x/\!\sigma)$
of the punctuation~$F_{\sigma,0}$ of~$F_\sigma$. 
In the finite case, since $[x f_\sigma(x)]'/\!f_\sigma(x)$ in condition~iii)
of Theorem~\ref{Thmscal} is just $\Lambda(x/\!\sigma)$, this theorem yields
  $ {\cal I}_{\rm s1}(F_\sigma) = \int \Lambda^2(x/\!\sigma)\,F_{\sigma,0}(dx) $,
  which is $ \int \Lambda^2(x)\,F_{0}(dx)= {\cal I}_{\rm s1}(F_\sigma)$;
that is, \eqref{Is1invariant} again. Therefore, in the finite case,
\begin{equation}
  {\cal I}_{\rm s}(F_\sigma)=\int \Lambda_\sigma^2 \,dF_{\sigma,0}\,,
   \quad 0<\sigma<\infty\,.
\end{equation}
the representation of~${\cal I}_{\rm s}(F_\sigma)$
in terms of the usual score function~$\Lambda_\sigma$,
\begin{equation}\label{scaMod}
      \Lambda_{\sigma}(x):=\frac{1}{\sigma}\,\Lambda\bigl(\frac{x}{\sigma}\bigr)=
      \frac{\partial}{\partial \sigma} \log f_\sigma(x) = -\frac{1}{\sigma} \Bigl(
      \!1 + \frac{x}{\sigma}\frac{f' ({x\over \sigma})}{f({x\over \sigma})}\,\Bigr) \:.
\end{equation}
\par
As an analogue to a lemma due to \citet{Ha:1972} in the location case,
\citet[Ch.2, Sec.3]{Sw:80} for an absolutely continuous~$F$ has shown that
conditions i)--iii) of Theorem~\ref{Thmscal} 
even imply \textit{$L_2$-differentiability} \citep[Def.~2.3.6]{Ri:94}
of the scale model,
\begin{equation} \label{l2def1}
  \bigl\Vert
  \sqrt{dF_{\sigma+t}}-\sqrt{dF_{\sigma}}(1+\Tfrac{1}{2}
   t\Lambda_{\sigma}) \bigr\Vert =\Lo(t) \qquad \mbox{as $t\to0$}
\end{equation}
at $\sigma=1$ and, by invariance, at any $0<\sigma<\infty$.
By definition, $L_2$-differentiability already entails that
  $\int \Lambda_\sigma^2 \,dF_\sigma<\infty$.
Setting $\Lambda(0):=0$, we may extend his result to $F(\{0\})>0$.
\begin{Prop}\label{Swens}
Assume that ${\cal I}_{\rm s1}(F)<\infty$\,.
Then the scale model $(F_\sigma)_{0<\sigma<\infty}$ is $L_2$-differentiable with
derivative~$\Lambda_\sigma$ at every $0<\sigma<\infty$\,.
\end{Prop}
\par $L_2$-differentiability of a parametric model implies an expansion
of the log-likelihhods, see e.g.\ \citet[Thm.~2.3.5]{Ri:94};
in our case, for each $h\in\R$,
\begin{equation} \label{LANiid}
\log dF_{\sigma+h\!/\!\sqrt{n}}^n/dF_{\sigma}^n = %
\Tfrac{1}{\sqrt{n}} \sum_{i=1}^n h^{\tau}\Lambda_{\sigma}(x_i)-\Tfrac{1}{2}
h^{\tau}{\cal I}_{{\rm s}}(F_\sigma)h+\Lo_{F_{\sigma}^n}(n^0)\,;
\end{equation}
that is, the scale model is \textit{locally asymptotically normal} (LAN).
LAN is the basis of asymptotic optimality results as H\'ajek's Asymptotic Convolution
Theorem and the Local Asymptotic Minimax Theorem, see e.g.\
\citet[Thm.'s~3.2.3, 3.3.8]{Ri:94} and
\citet[Thm.'s~8.8, 8.11]{VdW:98}.
\citet[17.3 Prop.~2]{LC:86} even shows that, in the i.i.d.\ setup, LAN
is equivalent to $L_2$-differen\-tiabi\-lity.
Thus we obtain the following result.
\begin{Prop}\label{equivprop}
The following statements are equivalent:
\begin{enumerate}
\item[i)] ${\cal I}_{\rm s}(F_\sigma)<\infty\,$ at some $0<\sigma<\infty$.
\item[ii)] The scale model is $L_2$-differentiable at some $0<\sigma<\infty$.
\item[iii)] The scale model has the LAN property \eqref{LANiid} at some $0<\sigma<\infty$.
\end{enumerate}
\end{Prop}
By invariance, the validity of each statement at one $\sigma$ implies its validity
at any other $0<\sigma<\infty$.
\appendix
\small
\makeatletter
\gdef\theThm{\@Alph\c@section.\arabic{Thm}}
\makeatother
\section{Proofs and Absolute Continuity}
\textit{Proof of Proposition~2.1}\enskip
The sup over a family of l.s.c., resp.\ convex, functions being l.s.c., resp.\ convex,
it suffices to show that, for each $\phi\in {\cal C}_{\rm c1}$, the reciprocal function
$V_1^{-1}(\phi,\,\cdot\,)$ from~\eqref{asyVarMsc}, is weakly l.s.c.\ and convex.
In this proof only, we pay a price for the simplifying convention $0/0:=0$.
\par
Let $F_n\to F$ weakly. Then $\int \phi^2\,dF_n\to \int \phi^2\,dF$.
First assume $\int \phi^2\,dF>0$. Then $\int \phi^2\,dF_n>0$ eventually,
and $V_1^{-1}(\phi,F_n)\to V_1^{-1}(\phi,F)$.
Secondly suppose that $\int \phi^2\,dF=0$. If also $\int x \,\phi'\,dF=0$, then
  $V_1^{-1}(\phi,F)=0 \le V_1^{-1}(\phi,F_n)$ for all~$n$.
If $\int x \,\phi'\,dF\ne 0$, then $\int \phi^2\,dF_n\to0$,
   $\int x \,\phi'\,dF_n\to \int x \,\phi'\,dF\ne0$, hence
   $V_1^{-1}(\phi,F_n)$ tends to $\infty =V_1^{-1}(\phi,F)$.
\par
Given $F_1$, $F_2$, $s\in (0,1)$, put $F=(1-s)F_1+sF_2$.
In case both $\int \phi^2\,dF_j>0$, we get
$ V_1^{-1}(\phi,F)\le
    (1-s)V_1^{-1}(\phi,F_1)+ s \,V_1^{-1}(\phi,F_2)$
    from \citet[Lemma~4.4]{Hu:81}.
Secondly, let $\int \phi^2\,dF_1=0<\int \phi^2\,dF_2$.
Then, if $\int x\,\phi'\,dF_1=0$, hence $V_1^{-1}(\phi,F_1)=0$,
and
  $V_1^{-1}(\phi,F)= s V_1^{-1}(\phi,F_2) = (1-s)0 + s \,V_1^{-1}(\phi,F_2)$.
If $\int x\,\phi'\,dF_1\ne0$, $V_1^{-1}(\phi,F_1)=\infty$ and
  $(1-s)\:\infty+s \, V_1^{-1}(\phi,F_2)\ge V_1^{-1}(\phi,F)$.
Thirdly, let both $\int \phi^2\,dF_j$ be zero.
Then, if also both $\int x\,\phi'\,dF_j=0$, we get $V_1^{-1}(\phi,F)=0$.
At least one $\int x\,\phi'\,dF_j$ nonzero implies that
  $(1-s)V_1^{-1}(\phi,F_1)+ s \,V_1^{-1}(\phi,F_2)=\infty$.
\hfill\qed
\begin{Lem}\sl \label{lemmadicht} For any finite measure~$F$ on~$\B$, the class
  ${\cal C}_{\rm c1}$ is dense in $L_2(F)$.  If $F(\{0\})=0$, the related class
     ${\cal D}_{\rm c1}:=\{\,x \mapsto x\,\phi'(x)\mid \phi\in {\cal C}_{\rm c1}\,\}$
is dense in $L_2(F)$.
There exist functions $0\le \phi_n\le1$ in~${\cal C}_{\rm c1}$ such that
  $\sup_{n,x}|x\,\phi_n'(x)|<\infty$, $\lim_n x\,\phi_n'(x)=0$,
  and $\phi_n(x) \uparrow1$, respectively $\phi_n(x) \downarrow 1_{\{x=0\}}$ pointwise.
\end{Lem}
\par \noindent \textit{Proof}\enskip
On the basis of Lusin's theorem, \citet[Thm.~3.14]{Ru:87}, it suffices to approximate
the indicator of bounded intervals~$(a,b]$. 
\par
For $\varepsilon \downarrow0$ one may choose functions
  $g_\varepsilon\in {\cal C}_{\rm c1}$ such that $0\le g_\varepsilon\le1$,
  $g_\varepsilon=1$ on $[a+\varepsilon,b]$,
  $g_\varepsilon=0$ on $(-\infty,a]\cup [b+\varepsilon,\infty)$.
Then $g_\varepsilon\to 1_{(a,b]}$ pointwise, and
  $g_\varepsilon\to 1_{(a,b]}$ in $L_2(F)$ by dominated convergence.
 \par
Concerning denseness of~${\cal D}_{\rm c1}$ in~$L_1(F_0)$, we may assume that $a>0$.
Drawing on the functions~$g_\varepsilon$
define $h_\varepsilon(x):=\int_{-\infty}^x y^{-1}g_\varepsilon(y)\,dy$.
Then $h_\varepsilon\in {\cal C}_{\rm c1}$
and,  as before,
  $x\,h_\varepsilon'=g_\varepsilon \to 1_{(a,b]}$ in $L_2(F_0)$.
\par
A possible choice of the functions~$\phi_n$, in the first case, is $\phi_n(x)=\phi(x/n)$,
based on the function $ 2 \,\phi(x) = 1 + \cos \bigl( ( |x|-\pi )_+ \land \pi\bigr)$, and,
in the second case, $\phi_n(x)=\phi(nx)$, where $2 \,\phi(x) = 1 + \cos ( |x|\land \pi)$.
\qed\smallskip \par
\textit{Absolute Continuity} \enskip 
From real analysis, e.g., \citet[Ch.8]{Ru:87}, we recall:
An $\R$-valued measure on the Borel $\sigma$-field~$\B$ of the real line is dominated by~$\lambda$,
the Lebesgue measure, iff its distribution function is absolutely continuous.
A function $f \colon\R\to\R$ is absolutely continuous, if for any $\varepsilon>0$ there
is a $\delta>0$ such that for any finite collection of disjoint segments $(a_i,b_i]$
of total length $\lambda \bigl( \bigcup (a_i,b_i] \bigr) <\delta$
it holds that $\sum_i|f(b_i)-f(a_i)|<\varepsilon$. 
Any absolutely continuous~$f$ has bounded variation on compact intervals $[a,b]$,
the derivative~$f'$ exists a.e.~$\lambda$,
and $f(b)-f(a)=\int_a^b f' \,d\lambda$ where $\int_a^b |f'|\,d\lambda<\infty$. 
Integrability $f'\in L_1(\lambda)$, implying bounded variation on~$\R$, and the
limit $f(a)\to0$ as $a\to -\infty$ require further conditions, respectively.
These are obviously satisfied in the location case for absolutely continuous
densities~$f$ such that ${\cal I}_{\rm l}(F)<\infty$ for $dF=f\,d\lambda$,
hence in particular $\int|f'|\,d\lambda<\infty$. 
If $f$ and~$g$ are absolutely continuous, so is their product $fg$
on any compact $[a,b]$.
Thus, integration by parts holds:
  $f(b)g(b)-f(a)g(a)= \int_a^b f'g\,d\lambda + \int_a^b fg'\,d\lambda$---a
  special case of \citet[Lemma~C.2.1]{Ri:94}.
\par \smallskip
\textit{Proof of Theorem 2.2}\enskip
First assume ${\cal I}_{\rm s1}(F)<\infty$. On~${\cal C}_{\rm c1}$ define
   $T(\phi):=-\int x\,\phi'\,dF$, which operator is well defined,
   because $\int \phi^2\,dF=0$, in view of Definition~\ref{FiSc},
   entails that $\int x\,\phi'\,dF=0$.
\par
Evaluated on~${\cal C}_{\rm c1}$, $T$ has operator norm~$\sqrt{\!{\cal I}_{\rm s1}(F)}\,$.
${\cal C}_{\rm c1}$ being dense in~$L_2(F)$, 
$T$ may be extended to~$L_2(F)$ keeping its norm.
By \textit{Riesz--Fr\'echet} there exists some $g\in L_2(F)$, whose norm equals the
operator norm of~$T$, such that
  $T(\phi)=\int \phi\,g\,dF$ for all $\phi\in L_2(F)$, hence
\begin{equation}\label{ixf=ifg}
            -\int x\,\phi'\,dF =
            \int \phi\,g\,dF \,,\qquad \phi\in {\cal C}_{\rm c1}\;.
\end{equation}
Inserting $\phi_n$ from Lemma~\ref{lemmadicht}, both choices, we obtain that,
in addition to $\int g^2\,dF={\cal I}_{\rm s1}(F)$,
\begin{equation}\label{Eigeng}
   \int g\:dF=0 \,,\qquad g(0)\,F(\{0\})=0
\end{equation}
In particular, the integrals in \eqref{ixf=ifg} and~\eqref{Eigeng} may be
restricted to $\R \setminus\{0\}$ . Define the function
\begin{equation}\label{deff}
            f(x):= \frac{1}{x}\int_{y\le x} g(y)\:F_0(dy) \,,\qquad x\ne0\;.
\end{equation}
Then, if $\phi_{-\infty}$ denotes the constant value of~$\phi\in {\cal C}_{\rm c1}$
left to the support of~$\phi'$,
   $\int \phi g\,dF = \int (\phi - \phi_{-\infty})\,g \,dF_0 $ and
        $\phi(x) - \phi_{-\infty}=\int_{0\ne y\le x} \phi'(y)\,\lambda(dy)$.
Due to compact support of~$\phi'$, and $g\in L_2(F_0)$, the product
  $ g(x)\,\phi'(y)$ is in $L_1(F_0(dx)\otimes \lambda(dy))$, and so
  $
   \int x\,\phi'\,dF_0 = -\iint_{x>y\ne0} g(x)\,\phi'(y)\,F_0(dx)\,\lambda(dy)
                       = \int yf(y) \,\phi'(y)\,\lambda(dy)
  $ by \textit{Fubini}; thus,
\begin{equation}\label{F0=fL}
   \int x\,\phi'(x)\,F_0(dx) = \int x\,\phi'(x)\,f(x)\,\lambda(dx)\,,
    \qquad \phi\in {\cal C}_{\rm c1}\;.
\end{equation}
By denseness of ${\cal D}_{\rm c1}$ in~$L_1(F_0)$, Lemma~\ref{lemmadicht},
the LHS determines~$F_0$.
As pointwise and dominated convergence
  $x\:h_\varepsilon'=g_\varepsilon \to 1_{(a,b]}$ has been established in that proof,
also $f\,d\lambda $ on the RHS is completely determined by~\eqref{F0=fL} if
  $f\,d\lambda$ is finite on any compact in $\R \setminus \{0\}$.
But
  $\int_A^B|f\,|\,d\lambda\le A^{-1}\int_A^B|xf(x)\,|\,\lambda(dx)$, which is
bounded by  $ (B/A-1) \int |g|\,dF_0<\infty$ for $A>0$, and likewise for $B<0$.
Thus we conclude from~\eqref{F0=fL} that
\begin{equation}
   dF_0=f\,d\lambda\;.
\end{equation}
Since $F_0$ is nonnegative, in fact $f\ge0$ a.e.~$\lambda$.
Absolute continuity of the function~$m$,
\begin{equation}
   m(x):= \int_{y\le x} g(y)\:F_0(dy) =
           \int_{y\le x} g(y)\,f(y\,) \,\lambda(dy)\;.
\end{equation}
follows from $\int |g|\,f\,d\lambda= \int |g|\,dF_0<\infty$.
As $m(x)=x\,f(x)$ for $x\ne0$, differentiability of~$f$ a.e.~$\lambda$ (for $x\ne0$)
is entailed by that of~$m$, and
\begin{equation}\label{identg}
   g(x)= 1 + x \:f'(x)\big/\!\!f(x) \,\qquad {\rm a.e.\ }F_0(dx)\;.
\end{equation}
This completes the identification of~$g$ under~$F$, and i)--iii) are proved.
\par
Conversely, assume i)--iii).
By~ii), $m(x)=x\,f(x)$ is absolutely continuous.
Differentiability of~$m$ at~$x\ne0$ implies that of~$f$, and $m'=f+xf'$.
For $\lambda$-densities, necessarily $\lambda(f=0, f'\ne0)=0$, hence also
    $\lambda(f=0,\: m'\ne0)=0$.
With $-\Lambda=m'\!\!/\!\!f = 1 + x\,f'\!\!/\!\!f\:$ a.e.~$F_0$,
     we have $\int|m'|\,d\lambda=\int|\,\Lambda|\,dF_0<\infty$ by~iii).
Thus,~$m$ and its measure $ m'\,d\lambda = -\Lambda\;dF_0$ are of bounded
variation on~$\R$.
\par
By H\"older inequality,
   $|m(y)-m(x)|^2\le |F(y)-F(x)|\int \Lambda^2\,dF_0$, so
   $m(x)$ for $x\to \infty$ is a Cauchy sequence.
But $\lim_{x\to \infty}m(x)$ must be zero since otherwise $f(x)\sim1/x$
for $x\to \infty$ would not integrate. 
The same holding for $x\to -\infty$,  
we obtain
\begin{equation}
       \int m'\,d\lambda=0.
\end{equation}
For $\phi\in {\cal C}_{\rm c1}$, the function $\phi - \phi_{-\infty}$ and
corresponding measure $ \phi'd\lambda$ have bounded variation on~$\R$.
Thus integration by parts in the general form of \citet[Lem.~C.2.1]{Ri:94}
yields $ \int \phi'm\,d\lambda = -\int \phi \,m'\,d\lambda$, such that
\begin{equation}
   \int x\,\phi'dF 
   = \int \phi'm\,d\lambda 
   = -\int \phi \,m'\,d\lambda = \int \phi\,\Lambda\,dF_0\;.
\end{equation}
Applying Cauchy-Schwarz, we get
\begin{equation}
  \Bigl(\:\int x\,\phi'dF \Bigr)^2 = \Bigl(\:\int \phi \,\Lambda\,dF_0 \Bigr)^2
  \le \int \phi^2\,dF_0 \int \Lambda^2\,dF_0\;,
\end{equation}
where $\int \Lambda^2\,dF_0$ is finite by~iii).
It follows that ${\cal I}_{\rm s1}(F)<\infty$. \hfill\qed
\par \medskip
%
\textit{Proof of Proposition~3.1}\enskip
We decompose
     $\|\sqrt{dF_{\sigma+t}}-\sqrt{dF_{\sigma}}(1+\Tfrac{1}{2}
       t\Lambda_{\sigma})\|$ into the following sum,
\begin{equation}
    \big\|\big(\sqrt{dF_{\sigma+t}}-\sqrt{dF_{\sigma}}(1+\Tfrac{1}{2}
    t\Lambda_{\sigma})\big)1_{\{0\}^c}\big\|
    + \big\|\big(\sqrt{dF_{\sigma+t}}-\sqrt{d F_{\sigma}}(1+\Tfrac{1}{2}
    t\Lambda_{\sigma})\big)1_{\{0\}}\big\|\:,
\end{equation}
The first summand is $\Lo(t)$ by \citet{Sw:80}. The second is $0$,
since $F_\sigma(\{0\})=F(\{0\})$ and $\Lambda_\sigma(0)=0$. \hfill\qed
\section*{Acknowledgements}
We thank two referees for their helpful comments.
\ifx\StProb\undefined
\section*{References}
\fi
%

%
%
\end{document}